\documentclass{amsart}
\usepackage{amsmath}
\usepackage{xy}
\usepackage{mathrsfs,amssymb,latexsym,graphicx,stmaryrd,xypic,nicefrac,multirow}
\author {Eric Goubault, Emmanuel Haucourt, Sanjeevi Krishnan}
\title {Covering space theory for directed topology}

\newcommand{\Tau}{\mathrm{T}}

\newcommand{\OP}[1]{#1^{\mathrm{op}}}
\newcommand{\BOX}{\mathbb{\square}}
\newcommand{\CATS}{\mathscr{C}}

\newcommand{\GENERIC}{\mathscr{C}}

\newcommand{\MONOIDS}{\mathscr{M}}
\newcommand{\POSETS}{\mathscr{P}}
\newcommand{\POSPACES}{\mathscr{P}}
\newcommand{\PREORDEREDSPACES}{\mathscr{Q}}

\newcommand{\SETS}{\mathscr{S}}

\newcommand{\STREAMS}{\mathscr{U}}
\newcommand{\SPACES}{\mathscr{T}}

\newcommand{\half}{\nicefrac{1}{2}}

\newcommand{\ra}{\rightarrow}
\newcommand{\hra}{\hookrightarrow}

\newcommand{\id}{{\mathrm{id}}}

\newcommand{\graph}[1]{\mathrm{graph}(#1)}

\newcommand{\sk}{\mathrm{sk}}

\newtheorem{thm}{Theorem}[section]

\newtheorem{cor}[thm]{Corollary}
\newtheorem{lem}[thm]{Lemma}

\newtheorem{prop}[thm]{Proposition}

\theoremstyle{definition}
\newtheorem{defn}[thm]{Definition}
\newtheorem{eg}[thm]{Example}

\newcommand{\cocubical}{{\oblong}}
\newcommand{\shell}{{\partial\oblong}}
\newcommand{\direalize}[1]{\upharpoonleft\!\!#1\!\!\downharpoonright}
\newcommand{\supp}{\mathrm{supp}}

\newtheorem{thm:main}{Theorem 3.12}
\newtheorem{cor:main}{Corollary 3.13}
\newtheorem{prop:calculation}{Proposition 2.23}
\newtheorem{prop:covering}{Proposition 2.24}

\begin{document}

\maketitle
\begin{abstract}
	The state space of a machine admits the structure of time.
	For example, the geometric realization of a precubical set, a generalization of an unlabeled asynchronous transition system, admits a ``local preorder'' encoding control flow.
	In the case where time does not loop, the ``locally preordered'' state space splits into causally distinct components.
	The set of such components often gives a computable invariant of machine behavior.  
	In the general case, no such meaningful partition could exist.
	However, as we show in this note, the locally preordered geometric realization of a precubical set admits a ``locally monotone'' covering from a state space in which time does not loop.
	Thus we hope to extend geometric techniques in static program analysis to looping processes.
\end{abstract}

\section{Introduction}
\begin{figure}\label{fig:discrete.to.continuous}
	\begin{center}
		\begin{tabular}{ccc}
			\includegraphics[width=30mm,height=30mm]{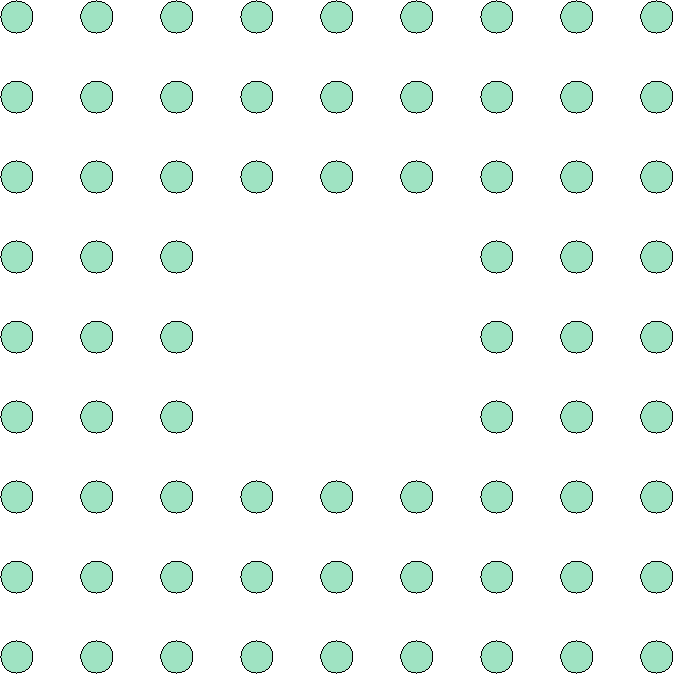}
			& \includegraphics[width=30mm,height=30mm]{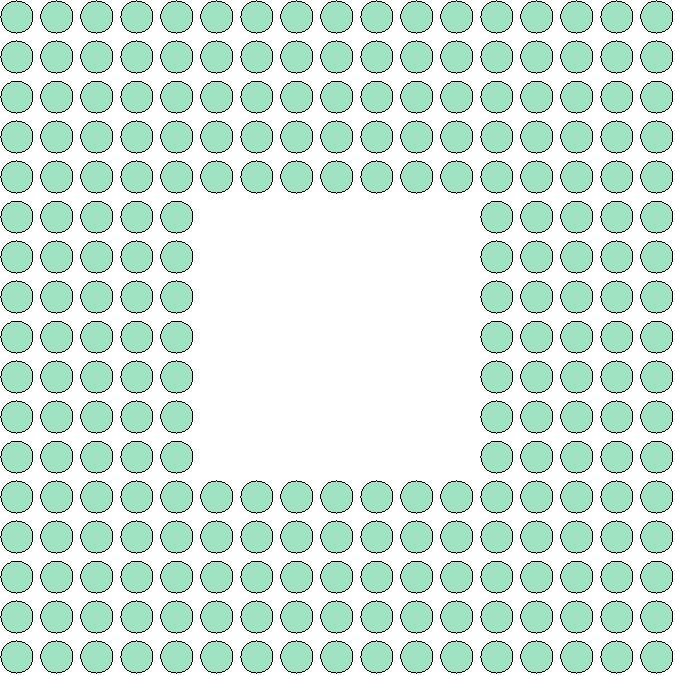}
			& \includegraphics[width=30mm,height=30mm]{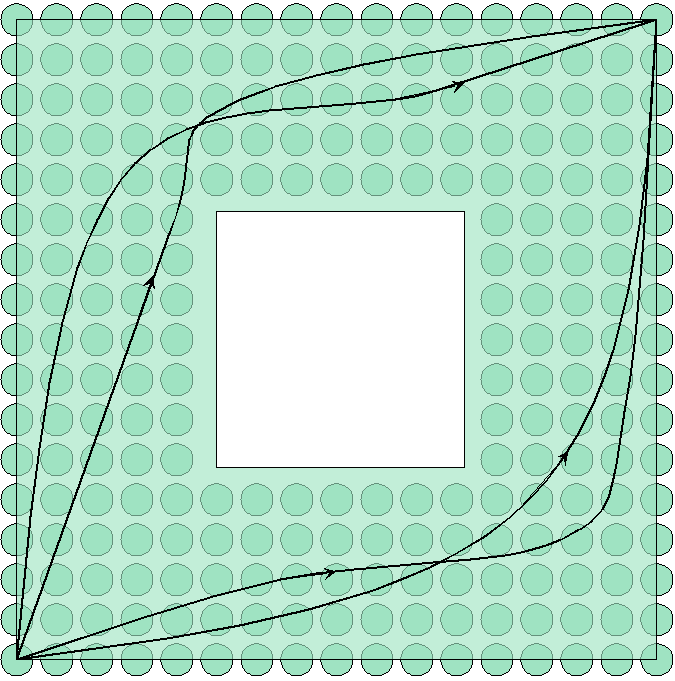}
		\end{tabular}
	\end{center}
	\caption{From discrete to continuous}
\end{figure}

The possible histories of a machine form a poset whose order encodes the progress of time.  
Meaningful machine behavior corresponds to properties of such posets independent of our measurement of time.
As the scales of time shrink, the sizes of posets explode. 
This combinatorial explosion often renders state space analyses incomplete in practice.
Figure \ref{fig:discrete.to.continuous} suggests that certain \textit{pospaces}, topological spaces equipped with suitable partial orders, represent the limits of such temporal refinement; chains of machine steps become monotone paths of evolutions.
Although these pospaces now contain an infinitude of histories, the \textit{directed homotopy theory} of \cite{fgr:ditop} on compact pospaces - more efficiently than order-theoretic analyses on finite but large posets - can extract the behavior of finite, terminating machines.
A static program analyzer, described in \cite{gh:alcool}, flags unsafe machine behavior by calculating \textit{component categories} of pospaces.

In the case where a finite machine does not terminate, its posets of histories are infinite and its associated pospace is non-compact.
Such partially ordered structures often reduce to a tractable size after the quotienting of histories ending at the same state.
We can recast locally finite posets as acyclic digraphs whose edges encode the partial orders, and we can enrich certain compact pospaces with the structure of ``local preorders.'' 
These new mathematical structures retain information about time flow after the quotienting process.
As the scales of time shrink, the sizes of digraphs explode.
``Locally preordered spaces'' represent the limits of such temporal refinement.
We follow \cite{krishnan:convenient} in calling these enriched spaces \textit{streams}.
Examples include \textit{stream realizations} $\direalize{X}$ of precubical sets $X$.
The directed homotopy theory of \cite{fgr:ditop}, adapted for compact streams, should circumvent the state space explosion problem in extracting the essential behavior of finite state machines.

\begin{figure}	
	\begin{center}
		\begin{tabular}{cc}
			\includegraphics[width=40mm,height=40mm]{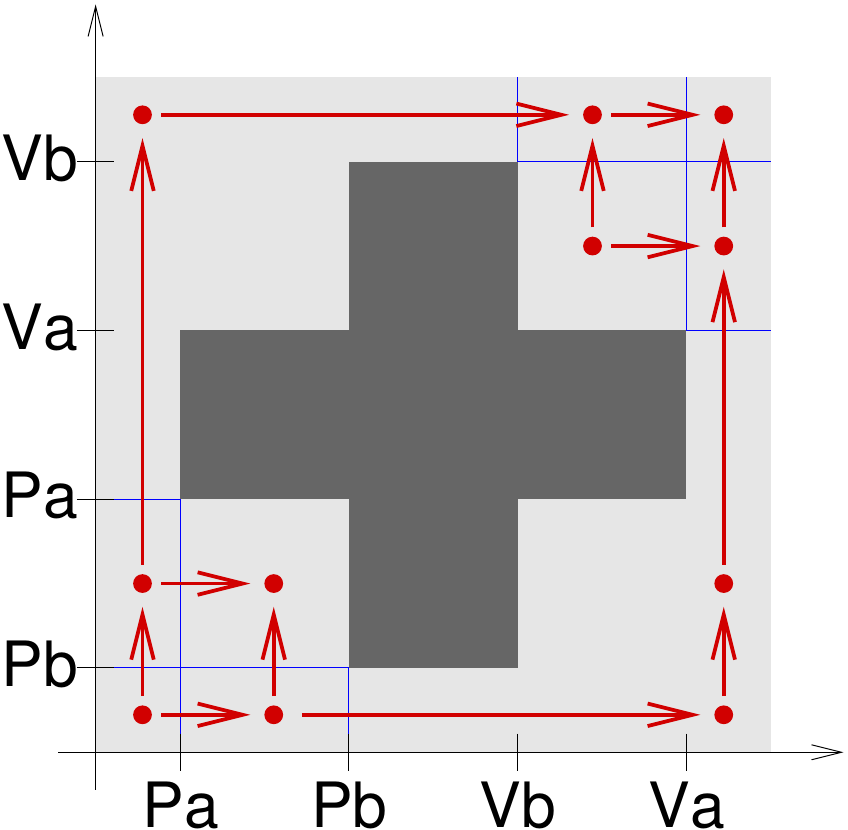} 
			& \includegraphics[width=40mm,height=40mm]{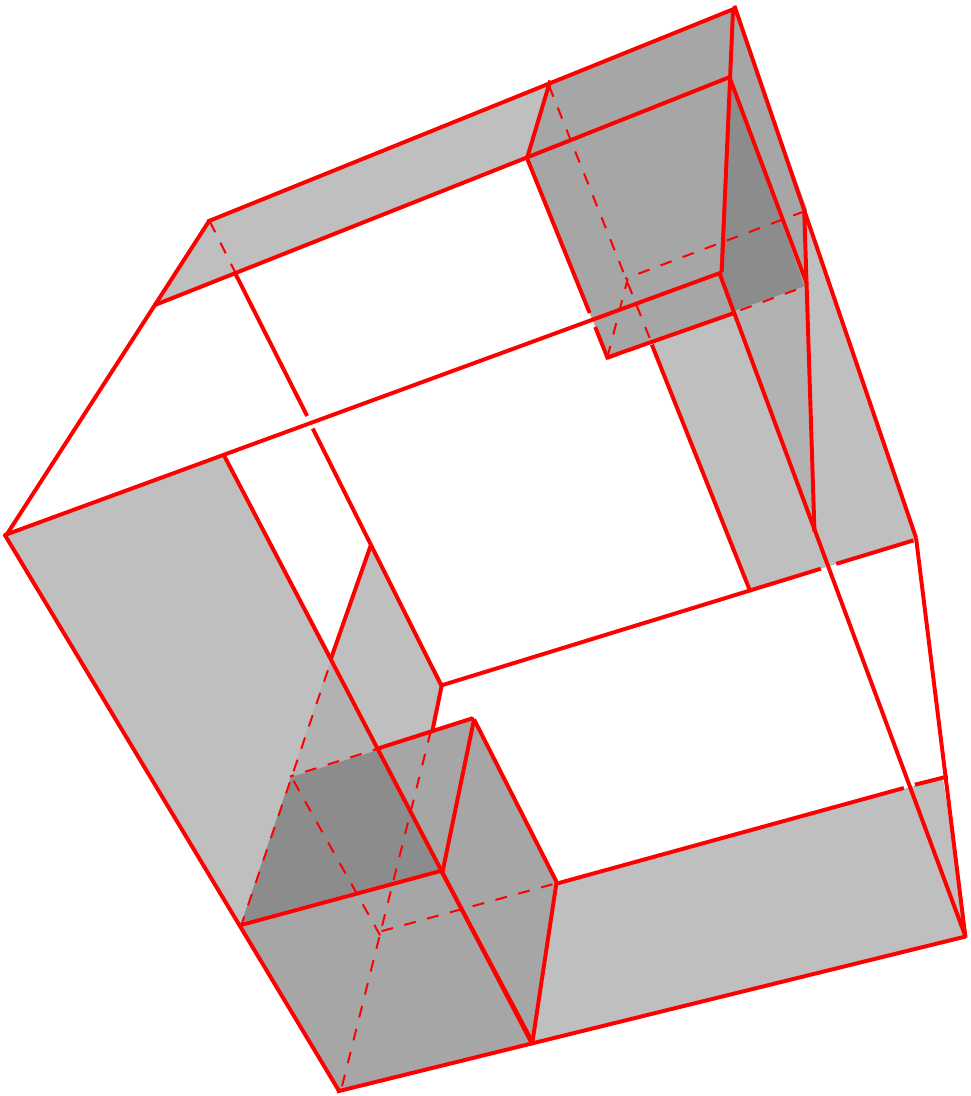}
			\vspace{.1in}\\
			(a) & (b)
		\end{tabular}
	\end{center}
	\caption{Tractable descriptions of machine behavior.
	The arrows in (a) and the lines in (b) comprise coherent choices of finitely many ``representative'' machine executions.}
	\label{fig:components}
\end{figure}

Path-components, the $0$-dimensional weak homotopy invariants of topological spaces, admit no straightforward and meaningful directed analogue.  
Fundamental groupoids of finite CW complexes admit finite skeletal subgroupoids whose objects correspond to path-components, while \textit{fundamental categories} of streams in nature rarely admit tractable skeletal subcategories, as observed in \cite{grandis:points}.
Yet, Figure \ref{fig:components} suggests that certain finite subcategories of fundamental categories of compact pospaces in nature appear to capture all possible subtleties in the dynamics of terminating machines.
Generalizations in \cite{fgrh:components, gh:components2} of isomorphisms in fundamental groupoids to ``causally insignificant'' morphisms in certain fundamental categories allow us to define tractable ``dipath-components'' of certain pospaces.

No such partitions of general streams could exist: a computer might reach the same state from different calculations.
In the case where the state stream of a machine arises as the quotient $X/G$ of a pospace $X$ by the transitive, free, and discrete action of a group $G$, the group $G$ and the ``dipath-components'' of $X$ together suffice in constructing meaningful geometric invariants on $X/G$.
In this note, our main result implies that stream realizations of connected precubical sets always assume such a form $X/G$.

\begin{thm:main}
	For each precubical set $B$, there exists a stream covering
	$$E\ra\;\direalize{B}$$
	such that the preorder $\leqslant_E$ is antisymmetric.
\end{thm:main}

For example, the quotient of the linearly ordered topological space $\vec{\mathbb{R}}$ of reals by the additive action of the group $\mathbb{Z}$ of integers yields the state stream $\vec{\mathbb{S}}$ of a cyclical process.
We call stream-theoretic analogues of the ``locally partially ordered spaces'' in \cite{fgr:ditop} \textit{vortex-free}.

\begin{cor:main}
	Stream realizations of precubical sets are vortex-free.
\end{cor:main}

We start in \S\ref{sec:ditop} by reviewing basic directed topology.  
In \S\ref{subsec:pospaces}, we recall definitions and examples of preordered spaces.
In \S\ref{subsec:streams}, we recall from \cite{krishnan:convenient} a convenient category of \textit{streams}, spaces equipped with preorderings $V\mapsto\leqslant_V$ of their open subsets satisfying properties reminiscent of the sheaf condition.
The forgetful functor $U:\STREAMS\ra\SPACES$ to the category of spaces is topological.
In particular, $\STREAMS$ is complete and cocomplete.
Moreover, compact pospaces whose maximal chains are connected admit unique compatible structures of streams.
In \S\ref{subsec:fundamental.categories}, we construct \textit{fundamental categories} of streams, in the process constructing \textit{fundamental monoids} $\tau_1(X,x)$ of based streams $(X,x)$.
We adapt covering space theory to the world of streams in \S\ref{subsec:coverings}.
A local preorder on a space $B$ pulls back along a covering map $E\ra B$ to define a \textit{stream covering}.
The usual path-lifting properties for coverings imply \textit{dipath}-lifting properties for stream coverings because $U$ is topological.  
In particular, our stream coverings are ``dicoverings'' in the sense of \cite{fajstrup:dicovering}.
A calculation follows.

\begin{prop:calculation}
	For each $x\in\mathbb{S}$, $\deg:\tau_1(\vec{\mathbb{S}},x)\cong\mathbb{N}$.
\end{prop:calculation}

Under mild hypotheses, we identify necessary and sufficient conditions for a universal cover of a stream to admit no loops in time.

\begin{prop:covering}
	The following are equivalent for a universal stream covering
	$$E\ra B$$
	over a path-ordered stream $B$.
	\begin{enumerate}
		\item The preorder $\leqslant_E$ is antisymmetric.
		\item\label{item:di.universal} Every stream map $\vec{\mathbb{S}}\ra B$ null-homotopic as a continuous function is constant.
	\end{enumerate}
\end{prop:covering}

In \S\ref{sec:precubical}, we turn our attention to the state streams of concurrent machines.
We review the basic definitions of precubical sets in \S\ref{subsec:precubical.sets}, described in \cite{goubault:transition} as generalizations of unlabeled, asynchronous transition systems.
Consider a precubical set $X$.
We enrich the geometric realization $|X|$, defined in \S\ref{subsec:realizations} as a certain CW complex, with the structure of a stream $\direalize{X}$ in \S\ref{subsec:stream.realizations}. 
Non-constant cellular dipaths on $\direalize{X}$ represent non-trivial cellular homology classes in $H_1(|X|;\mathbb{Z})$, by \cite[Example 2.8]{fajstrup:dicovering} and Lemma \ref{lem:dihomology}.
Non-constant dipaths on $\direalize{X}$ admit non-constant cellular approximations by Lemma \ref{lem:approx}.
The stream $B=\;\direalize{X}$ thus satisfies condition (\ref{item:di.universal}) of Proposition \ref{prop:covering}.
The main result ensues.

\section{Directed topology}\label{sec:ditop}

We review basic definitions of preordered spaces in \S\ref{subsec:pospaces}. 
We review basic definitions and properties of the category $\STREAMS$ of streams in \S\ref{subsec:streams}.
We construct fundamental categories and fundamental monoids in \S\ref{subsec:fundamental.categories}.
Lastly, we adapt covering space theory for streams in \S\ref{subsec:coverings}.

\subsection{Preordered spaces}\label{subsec:pospaces}
A \textit{preordered space} is a space $X$ equipped with a preorder $\leqslant_X$.on its points.
A \textit{pospace} is a preordered space $X$ whose preorder is a partial order whose graph $\graph{\leqslant_X}$ is closed in the standard product topology of $X\times X$.
In other words, a pospace is a poset topologized so that topological limits commute with inequalities.
Examples include the the real number line $\vec{\mathbb{R}}$ equipped with its standard total order.
Every pospace is Hausdorff by \cite[Proposition VI-1.4]{scott:lattices}.  

\begin{eg}
	A Hausdorff space with its trivial partial order is a pospace.
\end{eg}

Fix a preordered space $X$.  
A preordered space $A$ is a \textit{sub-preordered space} of $X$ if $A$ is a subspace of $X$ and $\graph{\leqslant_A}=\graph{\leqslant_X}\cap(A\times A)$.
Examples include the sub-preordered space $\vec{\oblong}[1]$ of $\vec{\mathbb{R}}$ whose underlying space is the unit interval $\mathbb{I}$.

A \textit{monotone map} is a (weakly) monotone and continuous function between preordered spaces.
Let $\PREORDEREDSPACES$ be the category of preordered spaces and monotone maps.
Colimits and limits in $\PREORDEREDSPACES$ are created by forgetful functors to suitable categories of spaces and proerdered sets.

\subsection{Streams}\label{subsec:streams}

We define \textit{streams} to be spaces equipped with coherent preorderings of their open subsets.
We summarize the basic definitions and properties of streams in this section, referring the reader to \cite{krishnan:convenient} for comparisons with the local pospaces of \cite{fgr:ditop} and the d-spaces of \cite{grandis:d}.

\begin{defn}
	A \textit{circulation} $\leqslant$ on a space $X$ is a function assigning to each open subset $V\subset X$ a preorder $\leqslant_V$ on $V$ such that for every collection $\mathcal{O}$ of open subsets of $X$, $\leqslant_{\bigcup\mathcal{O}}$ is the preorder on $\bigcup\mathcal{O}$ with smallest graph containing 
	$$\bigcup_{V\in\mathcal{O}}\graph{\leqslant_V}.$$
	A \textit{stream} $X$ is a space equipped with a circulation on it, which we always write as $\leqslant$.  
\end{defn}

\begin{eg}\label{eg:complex}
	For each open subset $V$ of the complex plane $\mathbb{C}$ and for all $z,z'\in\mathbb{C}$, define $z\leqslant_V z'$ if there exist non-decreasing paths $\alpha,\beta:\mathbb{I}\ra\mathbb{R}$ such that $\alpha e^{i\beta}$ defines a path $z\leadsto z'$ in $V$.
	The function $V\mapsto\leqslant_V$ turns $\mathbb{C}$ into a stream $\vec{\mathbb{C}}$.
\end{eg}

\begin{defn}
	Consider streams $X$ and $Y$.
	A \textit{stream map} is a continuous function 
	$$f:X\rightarrow Y$$
	satisfying $f(x)\leqslant_Vf(y)$ whenever $x\leqslant_{f^{-1}V}y$, for each open subset $V\subset Y$.
\end{defn}

Let $\STREAMS$ be the category of streams and stream maps between them and $\SPACES$ be the category of spaces and continuous functions between them.
We refer the reader to \cite{borceux:2} for the basic definitions and properties of topological functors.

\begin{prop}[{\cite[Proposition 5.8]{krishnan:convenient}}]\label{prop:topological}
	The forgetful functor 
	$$U:\STREAMS\rightarrow\SPACES$$
	is topological.
\end{prop}

In particular, $U$ creates limits and colimits by \cite[Proposition 7.3.8]{borceux:2}, and $\STREAMS$ is hence complete and cocomplete.

\begin{defn}
	Fix a stream $X$.
	A stream $A$ is a \textit{substream} of $X$ if inclusion defines a stream map $\iota:A\hookrightarrow X$ and each stream map $f$ whose image lies in $A$ corestricts to a dotted stream map making the following diagram commute.
	$$\xymatrix{
	A\ar[r]^{\iota}
	& X\\
	B\ar@{.>}[u]\ar[ur]_f
	}$$
	A \textit{stream embedding} $f:A\ra X$ is a stream map which corestricts to an isomorphism between $A$ and a substream of $X$.
\end{defn}

Fix a stream $X$.
Every subset of $X$ admits the unique structure of a substream of $X$, by Proposition \ref{prop:topological}.  
Generally, substreams of $X$ are difficult to describe explicitly.
However, open substreams of $X$ are open subspaces equipped with suitable restrictions of the circulation on $X$.

There exists a concrete forgetful functor $Q:\STREAMS\rightarrow\PREORDEREDSPACES$ sending each stream $X$ to its underlying space preordered by $\leqslant_X$.
The following observation allows us to calculate colimits and finite products of streams as colimits and finite products of underlying spaces, equipped with ``point-wise'' colimits and finite products of preordered sets.

\begin{prop}[{\cite[Lemma 5.5, Proposition 5.11]{krishnan:convenient}}]\label{prop:almost.topological}
	The functor 
	$$Q:\STREAMS\rightarrow\PREORDEREDSPACES$$
	preserves colimits and finite products.
\end{prop}

The direct verification of the axioms for streams and stream maps can be tedious.  
Let $\POSPACES$ be the category of pospaces whose maximal chains are connected, and all (weakly) monotone, continuous functions between them.
The following observation, a consequence of \cite[Propositions 4.7, 5.4, 5.11]{krishnan:convenient} and \cite[Propositions 1, 2, Theorem 5]{nachbin:order}, 
allows us to henceforth treat $\vec{\mathbb{R}}$ as a stream and connected sub-pospaces of $\vec{\mathbb{R}}$ as substreams.

\begin{thm}\label{thm:embed}
	There exists a full, concrete, and product-preserving embedding
	\begin{equation}\label{eqn:embed}
		\POSPACES\hookrightarrow\STREAMS
	\end{equation}
	sending each pospace to a unique stream sharing the same underlying space and underlying preordered set.
\end{thm}

Just as locally path-connected spaces are those spaces on which paths detect connectedness, \textit{path-ordered} streams are those streams on which \textit{dipaths} detect order-theoretic information.
Path-ordered streams essentially are the d-spaces of \cite{grandis:d} satisfying certain extra properties.

\begin{defn}
	Fix a stream $X$.
	A \textit{dipath} on $X$ is a stream map
	$$\gamma:\vec{\oblong}[1]\ra X,$$
	written $\gamma:x\leadsto y$ if $x=\gamma(0)$ and $y=\gamma(1)$.
	A stream is \textit{path-ordered} if there is a dipath $\gamma:x\leadsto y$ on $X$ whose image lies in $V$ whenever $x\leqslant_Vy$, for each open subset $V\subset X$ and each pair $x,y\in V$.
\end{defn}

\begin{eg}
	The stream $\vec{\mathbb{C}}$ from Example \ref{eg:complex} is path-ordered.
\end{eg}

Colimits of path-ordered streams are path-ordered, because dipaths concatenate.
Following \cite{grandis:d}, we define a \textit{vortex} to be the point of a stream having no partially ordered neighborhood.
\textit{Vortex-free} streams are analogues of the local pospaces of \cite{fgr:ditop}.

\begin{defn}
	Fix a stream $X$.
	A \textit{vortex} in $X$ is a point $x\in X$ such that
	$$\leqslant_V$$
	is not antisymmetric for every open neighborhood $V$ of $x$ in $X$.
	Call $X$ \textit{vortex-free} if it has no vortices.
\end{defn}

\begin{eg}
	In the stream $\vec{\mathbb{C}}$ from Example \ref{eg:complex}, $0$ is the unique vortex.
\end{eg}

An example of a vortex-free stream is $\vec{\mathbb{S}}$, the unit circle $\mathbb{S}=\{z\in\mathbb{C}\;|\;|z|=1\}$ equipped with the circulation $\leqslant$ defined by $z\leqslant_V z'$ if there exist reals $\theta\leq\theta'$ such that $z=e^{i\theta}$, $z'=e^{i\theta'}$, and $e^{i[\theta,\theta']}\subset V$.
The argument given in \cite[Example 2.8]{fajstrup:dicovering} adapts.

\begin{lem}[{\cite[Example 2.8]{fajstrup:dicovering}}]\label{lem:fajstrup}
	Fix a vortex-free stream $X$.  
	A stream map
	$$\vec{\mathbb{S}}\ra X$$
	is constant if and only if it is homotopic through stream maps to a constant map.
\end{lem}

\subsection{Fundamental categories}\label{subsec:fundamental.categories}
We adapt the construction in \cite{fgr:ditop} of fundamental categories and fundamental monoids, generalizing fundamental groupoids and fundamental groups.  

\begin{defn}
	Fix a stream $X$.  
	Let $\Tau_1X$ be the category having the points of $X$ as objects, equivalence classes $[\alpha]$ of dipaths homotopic relative $\{0,1\}$ through dipaths as the morphisms, and the functions
	$$s:[\alpha]\mapsto\alpha(0),\quad t:[\alpha]\mapsto\alpha(1),\quad\id:x\mapsto[t\mapsto x],\quad[\beta]\circ[\alpha]=[\alpha*\beta],$$
	as source $s$, target $t$, identity $\id$, and composition $\circ$ functions, respectively, where $\alpha*\beta$ denotes a concatenation of $\alpha$ with $\beta$.
	We call $\Tau_1X$ the \textit{fundamental category} of $X$.
\end{defn}

The construction $\Tau_1$ naturally extends to a functor $\Tau_1:\STREAMS\ra\CATS$ to the category $\CATS$ of small categories and functors.
For each stream $X$ and each $x\in X$, let $\tau_1(X,x)$ be the endomorphism monoid $(\Tau_1X)(x,x)$.
Recall that a monoid $M$ is \textit{inverse-free} if its only element admitting a left or right inverse is the unit.

\begin{lem}\label{lem:inverse-free}
	For each vortex-free stream $X$ and each $x\in X$, the monoid
	$$\tau_1(X,x)$$
	is inverse-free.
\end{lem}
\begin{proof}
	A representative dipath $\alpha$ of an element of $\tau_1(X,x)$ admitting a (left or right) inverse is constant by Lemma \ref{lem:fajstrup} because $\alpha:x\leadsto x$ concatenates (on the left or right) with a dipath $\beta:x\leadsto x$ on $X$ to form a dipath homotopic through dipaths to a constant map at $x$.
\end{proof}
 
The lemma generalizes to arbitrary streams $X$ and points $x$ which are not vortices.

Let $\STREAMS_\star$ be the category of based streams and stream maps preserving distinguished points and $\SPACES_\star$ be the category of based spaces and based continuous functions.
The construction $\tau_1$ naturally extends to a functor $\tau_1:\STREAMS_\star\ra\MONOIDS$ to the category $\MONOIDS$ of monoids and monoid homomorphisms.
The forgetful functor $U$ induces a natural transformation $U_*:\tau_1\ra\pi_1\circ U$.
Let $\deg$ denote both the standard degree map and the dotted homomorphism defined by the commutative diagram
$$\xymatrix@C=3pc{
\tau_1(\vec{\mathbb{S}},1)\ar[d]_{U_*}\ar@{.>}[r]^{\deg}
& \mathbb{N}\ar[d]^i\\
\pi_1(\mathbb{S},1)\ar[r]_{\deg}
& \mathbb{Z},
}$$
where the right vertical arrow denotes inclusion from the monoid $\mathbb{N}$ of natural numbers.

\subsection{Coverings}\label{subsec:coverings}
Covering space theory straightforwardly adapts.

\begin{defn}
	A \textit{stream covering} is a surjective stream map
	$$\rho:E\ra B$$
	such that $B$ admits a cover of open substreams whose preimages under $\rho$ are disjoint unions of open components on which $\rho$ restricts to define stream embeddings.
\end{defn}

A stream covering is a covering of underlying spaces.  
Conversely, a circulation on a space $B$ pulls back along a covering $\rho:E\ra B$ of spaces to define a circulation on $E$ turning $\rho$ into a stream covering by Proposition \ref{prop:topological}.
The following statements thus follow from their classical counterparts.

\begin{lem}\label{lem:path.lifting}
	Fix a stream covering $\rho:E\ra B$.
	\begin{enumerate}
		\item The stream $E$ is path-ordered if and only if $B$ is path-ordered.
		\item For each $e\in E$, every homotopy through dipaths on $B$ starting (ending) at $\rho(e)$ lifts under $\rho$ to a unique homotopy through dipaths starting (ending) at $e$.
	\end{enumerate}
\end{lem}

\begin{eg}
	Let $\vec{\mathbb{C}}$ be the stream from Example \ref{eg:complex}.
	Let
	$$f:\vec{\mathbb{R}}^2\ra\vec{\mathbb{C}}$$
	be the stream map sending $(x,y)$ to $xe^{iy}$.  
	The map $f$ is not a stream covering because $Uf$ is not a covering of underlying spaces, even though dipaths on $\vec{\mathbb{C}}$ lift under $f$ to dipaths on $\vec{\mathbb{R}}^2$.
\end{eg}

Lemma \ref{lem:path.lifting} allows us to make a non-trivial calculation.

\begin{prop}\label{prop:calculation}
	For each $x\in\mathbb{S}$, $\deg:\tau_1(\vec{\mathbb{S}},x)\cong\mathbb{N}$.
\end{prop}
\begin{proof}
	We take $x=1$ without loss of generality.
	The function $t\mapsto e^{\pi it}$ defines a universal stream covering $\vec{\mathbb{R}}\ra\vec{\mathbb{S}}$.  
	Dipaths $\gamma$ starting and ending at $1$ lift under this covering to unique dipaths starting at $0$ and ending at $\deg(\gamma)$ by Lemma \ref{lem:path.lifting}.
	Each pair $\alpha,\beta$ of dipaths on $\vec{\mathbb{R}}$ starting at $0$ and ending at the same integer are homotopic relative $\{0,1\}$ through dipaths $h_t=(1-t)\alpha+t\beta$, and therefore
	\begin{equation}
	U_*:\tau_1(\vec{\mathbb{S}},1)\ra\pi_1(\mathbb{S},1)=\mathbb{Z}
		\label{eqn:injective}
	\end{equation}
	is injective.  
	The image of (\ref{eqn:injective}) is therefore an inverse-free sub-monoid of $\mathbb{Z}$ by Lemma \ref{lem:inverse-free}.  
	The image of (\ref{eqn:injective}) contains $1\in\mathbb{Z}$ because the identity on $\mathbb{S}$ defines a stream map $\vec{\mathbb{S}}\ra\vec{\mathbb{S}}$.
	The calculation follows from the last two observations.
\end{proof}

We call a stream covering $p$ \textit{universal} if the covering $Up$ is universal.

\begin{prop}\label{prop:covering}
	The following are equivalent for a universal stream covering
	$$E\ra B$$
	over a path-ordered stream $B$.
	\begin{enumerate}
		\item The preorder $\leqslant_E$ is antisymmetric.
		\item Every stream map $\vec{\mathbb{S}}\ra B$ null-homotopic as a continuous function is constant.
	\end{enumerate}
\end{prop}
\begin{proof}
	Consider a universal stream covering 
	$$\rho:E\ra B.$$ 

	Suppose $\leqslant_E$ is antisymmetric. 
	Consider a stream map $\alpha:\vec{\mathbb{S}}\ra B$ such that $U(\alpha)$ is null-homotopic to a constant map.
	Then $U(\alpha)$ lifts to a continuous function $\mathbb{S}\ra E$, which must define a stream map $\tilde\alpha:\vec{\mathbb{S}}\ra E$ by Lemma \ref{lem:path.lifting}.
	Then $\tilde\alpha$ must be constant because $\leqslant_E$ is antisymmetric.
	Thus we conclude $\alpha$ is constant.
	
	Now suppose instead that every stream map $\vec{\mathbb{S}}\ra B$ null-homotopic as a continuous function is constant.
	Consider a stream map $\beta:\vec{\mathbb{S}}\ra E$.
	The continuous function $U(\beta)$, and hence $U(\rho\circ\beta)$, is null-homotopic relative $\{0,1\}$ because $UE$ is simply connected.  
	Thus $\rho\circ\beta$, and hence $\beta$, is constant by assumption.
	The preorder $\leqslant_E$ is antisymmetric because $E$ is path-ordered by Lemma \ref{lem:path.lifting}.
\end{proof}

\section{Cubical models}\label{sec:precubical}

We review the definitions of precubical sets in \S\ref{subsec:precubical.sets},  investigate cellular chain complexes of geometric realizations in \S\ref{subsec:realizations}, and present our main results in \S\ref{subsec:stream.realizations}.
We adopt the integral notation $\int_{\GENERIC}^x F(x,x)$ for the coend of a functor $F:\OP\GENERIC\times\GENERIC\ra\mathscr{D}$ and $S\cdot X$ for an $S$-indexed coproduct of an object $X$.

\subsection{Precubical sets}\label{subsec:precubical.sets}
We often model combinatorial polytopes as presheaves over subcategories of the category $\POSETS$ of posets and (weakly) monotone functions.
For example, simplicial sets are presheaves over the full sub-category of $\POSETS$ containing the non-empty, finite ordinals $[n]=\{0<\cdots<n\}$.
Let $\langle\BOX,\oplus,[0]\rangle$ be the smallest sub-monoidal category of the Cartesian monoidal category $\langle\POSETS,\times,[0]\rangle$ containing the functions $\delta_{-},\delta_{+}:[0]\ra[1]$ sending $0$ to $0$ and $1$, respectively.
Basic assertions about $\BOX$-morphisms tend to have inductive proofs.
We leave the proof of the lemma below as an exercise.

\begin{lem}\label{lem:steps}
	For each pair $\varepsilon\leqslant_{[1]^n}\varepsilon'$ of elements, there exist $\BOX$-morphisms
	$$\phi_1,\ldots,\phi_n\in\BOX([1],[1]^n)$$
	such that $\phi_1(0)=\varepsilon$, $\phi_n(1)=\varepsilon'$, and $\phi_i(1)=\phi_{i+1}(0)$ for $i=1,2,\ldots,n-1$.
\end{lem}

Let $c\SETS$ be the functor category of \textit{precubical sets}, functors $\OP\BOX\ra\SETS$ to the category $\SETS$ of sets and functions; and \textit{precubical functions}, natural transformations between precubical sets.
Examples of precubical sets include the functors
$$\BOX[n]=\BOX(-,[1]^n):\OP\BOX\ra\SETS,\quad n=0,1,2,\ldots$$

Fix a precubical set $X$.  
For each natural number $n$, we write $X_n$ for $X([1]^n)$ and for each $i\in\{1,2,\ldots,n\}$, we write $d_{-i},d_{+i}:X_n\ra X_{n-1}$ for the functions to which $X$ sends the $\BOX$-morphisms $[1]^{i-1}\oplus\delta_-\oplus[1]^{n-i}$ and $[1]^{i-1}\oplus\delta_+\oplus[1]^{n-i}$, respectively.  
We write $\sigma_*$ for the image of an $n$-cube $\sigma\in X_n$ under the natural bijection
$$X_n\cong c\SETS(\BOX[n],X).$$

A precubical set $A$ is a \textit{sub-precubical set} of $X$ if object-wise inclusions define a precubical function $A\hra X$.
For each natural number $n$ and each $\sigma\in X_n$, let $\sk_n(X)$ and $\langle\sigma\rangle$ be the minimal sub-precubical sets of $X$ satisfying
$$(\sk_n(X))_n=X_n,\quad\langle\sigma\rangle_n=\{\sigma\}$$
and for later convenience, let $\sk_{-1}X=\varnothing$.

\subsection{Geometric realizations}\label{subsec:realizations}
We can interpret a precubical set as the data of a CW complex as follows.
For each natural number $n$, let $\cocubical[n]=\mathbb{I}^n$.
The assignment $[1]^n\mapsto\cocubical[n]$ extends to a functor $\cocubical:\BOX\ra\SPACES$ linearly extending each $\BOX$-morphism $[1]^m\ra[1]^n$ to a continuous function ${\oblong}[m]\ra{\oblong}[n]$.
Let $|-|$ be the functor
$$|-|=\int_{\BOX}^{[1]^n}-_n\cdot\oblong[n]:c\SETS\ra\SPACES.$$

The functor $|-|$ preserves inclusions.  
We make the natural identification
$$(|\BOX[n]|,|\sk_{n-1}(\BOX[n])|)\cong(\oblong[n],\shell[n]),$$
where $\shell[n]$ is the subspace of $\oblong[n]$ consisting of all points having at least one coordinate equal to $0$ or $1$.

Consider a precubical set $X$.
The space $|X|$ is a CW complex, whose attaching maps of $n$-cells correspond to restrictions and corestrictions of $|\sigma_*|$ to functions $\partial\!\oblong[n]\ra|\sk_{n-1}(X)|$, for $n$-cubes $\sigma\in X_n$.
Each $x\in|X|$ inhabits the interior of a unique closed $n$-cell, the realization $|\langle\supp(x)\rangle|$ of a sub-precubical set of $X$ generated by a unique cube $\supp(x)$ of $X$.

We describe pieces of the cellular chain complex $C(|X|;\mathbb{Z})$, whose $n$-chains correspond to the free $\mathbb{Z}$-module $\mathbb{Z}[X_n]$ over the set $X_n$ of generators.
Consider $\sigma\in X_1$.  
The differential $\partial_1:\mathbb{Z}[X_1]\ra\mathbb{Z}[X_0]$, up to a sign change, satisfies 
$$\partial_1(\sigma)=(d_{+1}\sigma)-(d_{-1}\sigma),$$
We describe $\partial_2$ at the point-set level, following \cite{may:concise}.
Fix homeomorphisms
$$\varsigma_1:\partial\oblong[2]\cong\mathbb{S},\quad\varsigma_2:\oblong[1]/\shell[1]\cong\mathbb{S}$$
of various models of the circle by, say, the rules
$$\varsigma_1(x+\half,y+\half)=\frac{x+iy}{|x+iy|},\quad\varsigma_2(t)=e^{2\pi it}.$$
For each $\sigma\in X_1$, $|\sigma_*|$ passes to quotients to define an isomorphism
$$\oblong[1]/\shell[1]\ra|\langle\sigma\rangle|/|\sk_0\langle\sigma\rangle|$$
whose inverse extends to a continuous function 
$$\pi_X(\sigma):|\sk_1(X)|/|\sk_0(X)|\ra\oblong[1]/\shell[1]$$
sending all other points to the quotiented point.
Fix a $\theta\in X_2$.  
For each $\sigma\in X_1$, let $\lambda(\theta,\sigma)$ be the dotted continuous function making the diagram
$$\xymatrix@C=4pc{
\mathbb{S}\ar@{.>}[rrrr]^{\lambda(\theta,\sigma)}
& & & & \mathbb{S}\\
\partial\oblong[2]\ar[u]^{\varsigma_1}\ar[r]
& \oblong[2]\ar[r]_{|\theta_*|}
& |\sk_1(X)|\ar[r]
& |\sk_1(X)|/|\sk_0(X)|\ar[r]_{\pi_X({\sigma})}
& \oblong[1]/\partial\oblong[1]\ar[u]_{\varsigma_2}
}$$
commute, where the first and third bottom arrows denote inclusion and quotienting, respectively.
Straightforward calculations reveal that
$$
\deg(\lambda(\theta,d_{-1}\theta))=\deg(\lambda(\theta,d_{+2}\theta))=-1,\;\deg(\lambda(\theta,d_{+1}\theta))=\deg(\lambda(\theta,d_{-2}\theta))=+1.$$

The differential operator $\partial_2:\mathbb{Z}[X_2]\ra\mathbb{Z}[X_1]$, up to a sign change, satisfies
\begin{eqnarray*}
	\partial_2(\theta)
	&=& \sum_{\sigma\in X_1}\deg(\lambda(\theta,\sigma))(\sigma)\\
	&=& (d_{+1}\theta)+(d_{-2}\theta)-(d_{-1}\theta)-(d_{+2}\theta).
\end{eqnarray*}

Consequently, $1$-cycles homologous to $0$ are trivial if their coefficients have the same sign.

\begin{lem}\label{lem:dihomology}
	Fix a precubical set $X$.  
	The preimage of $0$ under the function
	$$Z_1(C(|X|;\mathbb{Z}))\cap\mathbb{N}[X_1]\ra H_1(C(|X|;\mathbb{Z}))$$
	sending a $1$-cycle to its homology class, where $Z_1(C(|X|;\mathbb{Z}))$ is the set of $1$-cycles of $C(|X|;\mathbb{Z})$ and $\mathbb{N}[X_1]$ is the set of $1$-chains having non-negative coefficients, is $0$.
\end{lem}

\subsection{Stream realizations}\label{subsec:stream.realizations}
Geometric realizations of precubical sets inherit structure induced from orientations of the $1$-cubes.
For each natural number $n$, let $\vec{\cocubical}[n]$ be the product stream ${\vec{\cocubical}[1])}^n$.
The assignment $[1]^n\mapsto\vec{\oblong}[n]$ uniquely extends to a dotted functor $\vec{\oblong}:\BOX\ra\SPACES$ making the following diagram commute.
$$\xymatrix{
& \STREAMS\ar[d]^{U}\\
\BOX\ar@{.>}[ur]^{\vec{\oblong}}\ar[r]_{\oblong}
& \SPACES
}$$

Thus the functor $\direalize{-}$, defined by
$$\direalize{-}=\int_{\BOX}^{[1]^n}-_n\cdot\vec{\oblong}[n]:c\SETS\ra\STREAMS,$$
lifts $|-|$ along $U$ in the commutative diagram
$$\xymatrix{
& \STREAMS\ar[d]^{U}\\
c\SETS\ar@{.>}[ur]^{\direalize{\;-\;}}\ar[r]_{|-|}
& \SPACES,
}$$
by Proposition \ref{prop:topological} and Theorem \ref{thm:embed}.

\begin{defn}
	The \textit{stream realization} of a precubical set $X$ is $\direalize{X}$.
\end{defn}

\begin{lem}\label{lem:path.ordered}
	Stream realizations of precubical sets are path-ordered.
\end{lem}
\begin{proof}
	Stream realizations are colimits of the path-ordered streams $\vec{\oblong}[n]$.
\end{proof}

The $1$-skeleton of a precubical set encodes the order-theoretic relationships between vertices in a stream realization.
We record a special case.  
For each partial order $\preccurlyeq_X$ on a set $X$ and each subset $A\subset X$, let $(\preccurlyeq_X)_{\restriction A}$ denote the partial order on $A$ having graph $\graph{\preccurlyeq_X}\cap(A\times A)$.  

\begin{lem}\label{lem:restrict}
	For each natural number $n$, 
	$$(\leqslant_{\direalize{\;\BOX[n]\;}})_{\restriction\;\direalize{\;\sk_0(\BOX[n])\;}}=(\leqslant_{\direalize{\;\sk_1(\BOX[n])\;}})_{\restriction\direalize{\;\sk_0(\BOX[n])\;}}.$$
\end{lem}
\begin{proof}
	The lemma follows from Proposition \ref{prop:almost.topological} and Lemma \ref{lem:steps} because
	$$\leqslant_{[1]^n}=(\leqslant_{\vec\oblong[n]})_{\restriction[1]^n}.$$
\end{proof}

Dipaths on stream realizations admit cellular approximations.
Our proofs resemble the arguments used in \cite{fajstrup:approx}.

\begin{lem}\label{lem:baby.approx}
	Fix an integer $n>0$ and a dipath $\alpha$ on $\direalize{\BOX[n]}$.
	There is a dipath
	$$\beta:\supp(\alpha(0))_*(0,\ldots,0)\leadsto\supp(\alpha(1))_*(0,\ldots,0)$$
	on $\direalize{\BOX[n]}$ such that $U(\beta)$ is cellular.
	Moreover, $\beta$ is non-constant if $\supp(\alpha(1))$ inhabits $\sk_{n-1}(\BOX[n])$ and $\langle\supp(\alpha(\nicefrac{1}{2}))\rangle=\BOX[n]$.
\end{lem}
\begin{proof}
	Let $\chi_S$ be the characteristic function of a subset $S\subset[0,1]$.
	Let $\pi_i:\oblong[n]\ra\oblong[1]$ be the $i$th projection function.
	For each $i$,
	\begin{eqnarray*}		
		\pi_i(\supp(\alpha(0))_*(0,\ldots,0))
		&=&
		\chi_{\{1\}}(\pi_i(\alpha(0)))\\
		&\leqslant_{[1]^n}& \chi_{\{1\}}(\pi_i(\alpha(1)))\\
		&=& \pi_i(\supp(\alpha(1))_*(0,\ldots,0)).
	\end{eqnarray*}
	
	Lemmas \ref{lem:path.ordered} and \ref{lem:restrict} imply that there exists a cellular path defining a dipath
	$$\beta:\supp(\alpha(0))_*(0,\ldots,0)\leadsto\supp(\alpha(1))_*(0,\ldots,0).$$
	
	Assume the hypothesis of the second sentence in the statement of the lemma.
	Then $\beta(0)=(0,\ldots,0)$ because for each $i$, $\pi_i(\alpha(0))\leq\pi_i(\alpha(\half))<1$.  
	For each $i$, $0<\pi_i(\alpha(\half))\leq\pi_i(\alpha(1))$.  
	There exists a $j$ such that $\pi_j(\alpha(1))\in\{0,1\}$, and hence $\pi_j(\alpha(1))=1$.
	Thus $\beta(1)\neq\beta(0)$.
\end{proof}

\begin{lem}\label{lem:approx}
	Fix a precubical set $X$ and a dipath $\alpha$ on $\direalize{X}$.
	There is a dipath
	$$\beta:\supp(\alpha(0))_*(0,\ldots,0)\leadsto\supp(\alpha(1))_*(0,\ldots,0)$$
	on $\direalize{X}$ such that $U(\beta)$ is cellular and $U(\alpha)\sim U(\beta)$.
	We can take $\beta$ to be non-constant if $\alpha$ is non-constant and $\alpha(0)=\alpha(1)$.
\end{lem}
\begin{proof}
	By reparametrization of $\alpha$, we assume that there exists an integer $k>0$ such that for each $i=1,\ldots,k$, there exist natural number $d_i$ and $\sigma_i\in X_{d_i}$, such that $\alpha(\nicefrac{i-1}{k},\nicefrac{i}{k})\subset|\langle\sigma_i\rangle|\setminus|\sk_{d_i-1}\langle\sigma\rangle|$.
	We take $k$ to be minimal without loss of generality.
	Let $\alpha_i$ be the restriction of $\alpha$ to the substream $[\nicefrac{i}{k},\nicefrac{i+1}{k}]$ of $\vec{\oblong}[1]$.  

	There exist unique continuous functions $\tilde{\alpha}_i$ making the diagrams
	\begin{equation*}
		\xymatrix{
		& \direalize{\BOX[d_i]}\ar[d]^{\direalize{\;(\sigma_i)_*\;}} \\
		[\nicefrac{i-1}{k},\nicefrac{i}{k}]\ar@{.>}^{\tilde\alpha_i}[ur]\ar[r]_{\quad \alpha_i}
		& \direalize{X}
		}
	\end{equation*}
	of underlying spaces commute.  
	The restrictions of the $\tilde{\alpha}_i$'s to the sub-preordered spaces $(\nicefrac{i-1}{k},\nicefrac{i}{k})$ of $Q[\nicefrac{i-1}{k},\nicefrac{i}{k}]$ are monotone because the corestrictions $\direalize{\BOX[d_i]}\ra\direalize{\langle\sigma_i\rangle}$ of the $\direalize{(\sigma_i)_*}$'s restrict and corestrict to isomorphisms $\direalize{\BOX[d_i]}\setminus{\direalize{\sk_{d_i}\BOX[d_i]}}\ra\direalize{\langle\sigma_i\rangle}\setminus\direalize{\sk_{d_i}\langle\sigma_i\rangle}$ of open substreams by Proposition \ref{prop:almost.topological}.
	The $\tilde{\alpha}_i$'s define monotone maps $Q[\nicefrac{i-1}{k},\nicefrac{i}{k}]\ra Q\direalize{\BOX[d_i]}$ because the $\tilde\alpha_i\times\tilde\alpha_i$'s send dense subsets of the compact spaces $\graph{\leqslant_{[\nicefrac{i-1}{k},\nicefrac{i}{k}]}}$ into the compact and Hausdorff spaces $\graph{\leqslant_{\direalize{\BOX[d_i]}}}$.
	Thus the $\tilde\alpha_i$'s in fact define stream maps by Theorem \ref{thm:embed}.

	There exist stream maps $\tilde\beta_i:[\nicefrac{i-1}{k},\nicefrac{i}{k}]\ra\direalize{\BOX[d_i]}$ such that
	$$\tilde\beta_i(\nicefrac{i-1}{k},\nicefrac{i}{k})\subset\direalize{\sk_1\BOX[d_i]},\quad\tilde\beta_i(\nicefrac{i+\varepsilon-1}{k})=\supp(\tilde\alpha_i(\nicefrac{i+\varepsilon-1}{k}))_*(0,\ldots,0),\;\;\varepsilon=0,1$$
	by Lemma \ref{lem:baby.approx}.
	Let $\tilde{h}_i:U(\tilde\alpha_i)\sim U(\tilde\beta_i)$ be the homotopies defined by
	$$\tilde{h}_i(x,t)=(1-t)\tilde\alpha_i(x)+t\tilde\beta_i(x).$$
	The composites $\beta_i=|(\sigma_i)_*|\circ\tilde\beta_i$ agree on their overlap because 
	\begin{eqnarray*}
		\beta_i(\nicefrac{i}{k}) 
		&=& |(\sigma_i)_*|(\supp(\tilde\alpha_i(\nicefrac{i}{k}))_*(0,\ldots,0))\\
		&=& \supp(\alpha_i(\nicefrac{i}{k}))_{*}(0,\ldots,0)\\ 
		&=& \supp(\alpha_{i+1}(\nicefrac{i}{k}))_{*}(0,\ldots,0)\\ 
		&=& |(\sigma_i)_*|(\supp(\tilde\alpha_{i+1}(\nicefrac{i}{k}))_*(0,\ldots,0))\\
		&=& \beta_{i+1}(\nicefrac{i}{k}).
	\end{eqnarray*}
	Thus the composites $h_i=|(\sigma_i)_*|\circ\tilde{h}_i$ agree on their overlap to define our desired homotopy $h$.
	
	Suppose $\alpha$ is non-constant and $\alpha(0)=\alpha(1)$.
	Then $\tilde\alpha_1(\nicefrac{1}{k})\in|\sk_{d_1-1}\BOX[d_1]|$ because: in the case $k=1$, $\alpha_1(0)=\alpha_1(1)$; in the case $k>1$, $\langle\sigma_1\rangle\cap\langle\sigma_2\rangle\subset\sk_{d_1-1}\langle\sigma_1\rangle$ by the minimality of $k$.
	And $\langle\supp(\tilde\alpha_1(\half))\rangle=\BOX[d_1]$ by our assumption on $\alpha$.
	Thus $\tilde\beta_1$, and hence $h(-,1)$, are non-constant by Lemma \ref{lem:baby.approx}.
\end{proof}

We can translate properties of cellular $1$-cycles of $|X|$ having positive coefficients into order-theoretic properties of $\direalize{X}$.

\begin{lem}\label{lem:null.homotopic}
	Consider a precubical set $X$.
	A stream map
	$$\gamma:\vec{\mathbb{S}}\ra\direalize{X},$$
	is constant if and only if $U(\gamma):\mathbb{S}\ra|X|$ is null-homotopic.
\end{lem}
\begin{proof}
	Suppose $U(\gamma)$ is non-constant.
	We give $\mathbb{S}$ the structure of a CW complex with one vertex at $1$, and we can take $U(\gamma)$ to be cellular by Lemma \ref{lem:approx}.
	For each $\sigma\in X_1$, let $\lambda'(\gamma,\sigma)$ be the dotted stream map making the diagram
	$$
		\xymatrix@C=4pc{
		\vec{\mathbb{S}}\ar@{.>}[rr]^{\lambda'(\gamma,\sigma)}\ar[d]_{\gamma}
		& & \vec{\mathbb{S}}\\
		\direalize{\sk_1(X)}\ar[r]
		& \direalize{\sk_1(X)}/|\sk_0(X)|\ar[r]_{\quad\vec{\pi}_X(\sigma)}
		& \vec{\oblong}[1]/\partial\oblong[1]\ar[u]_{\vec{\varsigma}_2}
}
		\label{eqn:dihomology.class}
	$$
	commute, where the bottom left arrow is a quotient stream map and $\vec{\varsigma}_2,\;\vec{\pi}_X(\sigma)$ are the stream maps defined by the respective functions $\varsigma_2,\;\pi_X(\sigma)$ of underlying spaces.
	For each $\sigma\in X_1$, $\lambda'(\gamma,\sigma)(1)=1$ because $\gamma$ is cellular and $\deg(\lambda'(\gamma,\sigma))$ is non-negative by Proposition \ref{prop:calculation}. 
	Let
	$$c(\gamma)=\sum_{\sigma\in X_1}\deg(\lambda'(\gamma,\sigma))(\sigma)\in\mathbb{Z}_1(C(|X|;\mathbb{Z}))\cap\mathbb{N}[X_1].$$
	There exists $\sigma_0$ such that $\deg(\lambda'(\gamma,\sigma_0))\neq 0$ by Lemma \ref{lem:fajstrup}, because $\gamma$ is non-constant.
	If $U(\gamma)$ were null-homotopic, then $c(\gamma)$ would be homologous to $0$ by the homotopy invariance of cellular homology, contradicting Lemma \ref{lem:dihomology}.
\end{proof}

The main result follows from Proposition \ref{prop:covering} and Lemmas \ref{lem:path.ordered}, \ref{lem:null.homotopic}.

\begin{thm}\label{thm:main}
	For each precubical set $B$, there exists a stream covering
	$$E\ra\;\direalize{B}$$
	such that the preorder $\leqslant_E$ is antisymmetric.
\end{thm}

\begin{cor}
	Stream realizations of precubical sets are vortex-free.
\end{cor}

\section{Conclusion}

We have thus observed that the connected state spaces of generalized asynchronous transition systems form the orbit streams of pospaces equipped with transitive, free, and discrete actions of groups.
Consider a connected such state stream $B$, choose a basepoint $b\in B$, and let $G=\pi_1(UB,b)$.
There exists a universal covering $E\ra B$ such that $\leqslant_E$ is antisymmetric by Theorem \ref{thm:main}.  
In later work, we hope to extract meaningful machine dynamics from $B$ as the $G$-orbits of a $G$-algebraic gadget $I(E)$, constructed on the $G$-pospace $E$.

\section{Acknowledgements}
The authors appreciate helpful discussions on the subject with Lisbeth Fajstrup and Martin Raussen.

\end{document}